\documentclass[12pt]{article}
\usepackage{amsmath,amsthm,epsfig,epsf,psfrag,natbib}
\usepackage{graphicx,amssymb}
\usepackage{enumerate}
\usepackage[latin1]{inputenc}

\usepackage{mathrsfs}
\usepackage{dcolumn}
\usepackage{amscd,amsfonts,amssymb}
\usepackage{amsopn}

\usepackage{amsbsy,lscape,epstopdf}
\usepackage{color,latexsym,amsfonts,bbm,amscd}

\usepackage{url} 

\newcommand{\blind}{1}

\newcommand{\be}{\begin{equation}}
\newcommand{\ee}{\end{equation}}
\newtheorem{theorem}{Theorem}
\newtheorem{lemma}{Lemma}

\newtheorem{remark}{Remark}

\newcommand{\ba}{\begin{eqnarray}}
\newcommand{\ea}{\end{eqnarray}}
\newcommand{\bas}{\begin{eqnarray*}}
\newcommand{\eas}{\end{eqnarray*}}

\newcommand{\uu}{{\bf{u}}}
\newcommand{\vv}{{\bf{v}}}

\newcommand{\yy}{{\bf{y}}}

\newcommand{\diff}{\mathrm{d}}

\def\tr {{\bf tr}}

\newcommand{\greekbold}[1]{\mbox{\boldmath $#1$}}

\newcommand{\bbeta}{\greekbold{\beta}}
\newcommand{\bbetahat}{\widehat{\greekbold{\beta}}}

\newcommand{\bmu}{\greekbold{\mu}}
\newcommand{\muhat}{\hat{\mu}}

\newcommand{\bmus}{\greekbold{\scriptstyle \mu}}

\newcommand{\bSigma}{\mbox{\boldmath $\Sigma$}}
\newcommand{\bSigmas}{\greekbold{\scriptstyle \Sigma}}

\newcommand{\btheta}{\mbox{\boldmath $\theta$}}

\newcommand{\pr}{{\mbox{{Pr}}}}

\newcommand{\bq}{\mbox{\bf q}}

\newcommand{\cS}{\mbox{$\mathcal S$}}

\newcommand{\bF}{\mbox{\bf F}}

\newcommand{\bX}{\mbox{\bf X}}
\newcommand{\bY}{\mbox{\bf Y}}

\newcommand{\bS}{\mbox{\bf S}}
\newcommand{\bI}{\mbox{\bf I}}

\newcommand{\ind}{\mathbbm{1}}

\begin{document}

\def\spacingset#1{\renewcommand{\baselinestretch}%
{#1}\small\normalsize} \spacingset{1}

\if1\blind
{
  \title{\bf Multi-parameter One-Sided Monitoring Test}
  \author{Guangyu Zhu\thanks{
    The authors gratefully acknowledge funding from NSERC Grant RGPIN-2014-03743,
    a Collaborative Research and Development Grant from NSERC and FPInnovations,
    and the ``a thousand talents'' program through Yunnan University.
    We are also indebted to the Forest Products Stochastic Modelling Group 
    centered at the University of British Columbia (UBC): members of this group 
    from FPInnovations in Vancouver,  Simon Fraser University, and UBC provided 
    stimulating discussions of the long-term monitoring program to which this paper contributes.}\\
    Department of Statistics, University of British Columbia, Canada\\
        Jiahua Chen\\
    Research Institute of Big Data, University of Yunnan, China\\
    Department of Statistics, University of British Columbia, Canada}
  \maketitle
} \fi

\if0\blind
{
  \bigskip
  \begin{center}
    {\LARGE\bf Multi-parameter One-Sided Monitoring Test}
\end{center}
  \medskip
} \fi

\begin{abstract}
Multi-parameter one-sided hypothesis test problems arise naturally in many applications.
We are particularly interested in effective tests
for monitoring multiple quality indices in forestry products. 
Our search reveals that there are many effective statistical methods 
in the literature for normal data, and that they can easily be adapted for non-normal data.
We find that the beautiful likelihood ratio test is unsatisfactory,
because in order to control the size, it
must cope with the least favorable distributions at the cost of power. 
In this paper, we find a novel way to slightly ease
the size control, obtaining a much more powerful test.
Simulation confirms that the new test retains good control of
the type I error and is markedly more powerful than the likelihood ratio test
as well as many competitors based on normal data.
The new method performs well in the context of monitoring multiple quality indices. 
\end{abstract}

\vfill

\newpage
\spacingset{1.45} 
\section{Introduction}
\label{sec:intro}

The research problem in this paper is motivated by an application.
The reliability of a wood structure heavily depends on the mechanical strength
of its component wood. 
It is important to closely monitor the dynamic wood strength distribution 
of solid lumber over time.
This is done through data collected via a random sample from the target populations
and the subsequent data analysis.
A few weak components have potentially severe consequences for the
structure, so the lower quantiles of the strength distribution have received the most attention.
See the lumber-quality monitoring procedures specified
in the American Society for Testing and Materials (ASTM) Standard D1990
\citep{astm2002standard}. 
This is also evident from the recent report by \cite{verrill2015simulation}, 
which examined the performance of various tests in the context of
5\% quantiles.

Clearly, even if the strength distribution of the wood product meets
the quality standard for the lower quantiles, the median or mean
strengths could be significantly lower than the norm. The reliability of the
structure could still be seriously compromised. 
This suggests the need to develop a monitoring test 
procedure for several quality indices simultaneously.
We aim to draw the attention of practitioners to this need and 
to develop an effective and easy-to-use test procedure.

The application easily translates into a statistical question.
We wish to statistically detect potential danger arising when the values of several 
user-selected parameters fall below well-established standards.
In other words, we seek a test for multi-parameter one-sided null
and alternative hypotheses. More abstractly, suppose we have a sample
from distribution $F$, and $\theta$ is a vector-valued parameter
or functional of $F$. We wish to test the hypothesis
\be
\label{H1}
H_0: \btheta \geq \btheta^*
\mbox{
against the alternative }
H_a: \btheta \not \geq \btheta^*
\ee
for a specific known vector $\btheta^*$,
where the inequality is interpreted to be component-wise.
Because of the invariance property, without loss of generality, we may
take $\btheta^* = {\bf 0}$; this will be assumed hereafter
unless otherwise indicated. 
The dimension of $\btheta$ will be denoted as $p$.
Clearly, many existing tests can easily be adapted to
this problem. 
However, we suggest that none of them seem to 
exactly fit, and additional research is needed. 


Under the normal model, the likelihood ratio test (LRT)
provides standard solutions to the current pair of opposing hypotheses 
$H_0$ and $H_a$ and similarly formulated pairs of opposing hypotheses.
Statisticians must determine the appropriate rejection region to ensure that the 
LRT has the size specified by the user. Along this line,
\cite{robertson1988order} worked out the solution to the 
LRT problem for the case where $\bSigma$ is known to be $\bI$. 
\cite{perlman1969one} solved the LRT problem where $\bSigma$ is unknown.

By the standard definition in mathematical statistics, the size of a test
is the supremum of its type I error. When the null hypothesis is composite,
i.e., it contains many distributions,
the size of the test is the type I error in the worst scenario, 
or at the least favorable
null distribution. Controlling the size of the test
can therefore lead to a pessimistic procedure: the type I error under the likely true 
data-generating distribution is far below the size of the test that leads to
compromised power. This is particularly true for the LRT for 
multi-parameter one-sided hypotheses. 
\cite{perlman2003validity} 
and \cite{perlman2006some} examined the
rejection region of the LRT in many situations and developed more
powerful tests accordingly. Such research is often motivated by
medical studies, where
the aim is often to assess whether a therapy has a beneficial effect on multiple outcomes 
simultaneously relative to a control. 
The specifics of these one-sided hypotheses vary depending on the medical problem.
For instance, \cite{o1984procedures} and \cite{tang1989approximate} 
proposed and extended a generalized least-squares test
that is most powerful when the true population mean is near
a specific line in the alternative space. 
In clinical studies with multiple outcomes, researchers may wish to
confirm that a new treatment is superior in at least one of
the outcomes and equivalent on the rest of the outcomes,
in comparison with the control.
\cite{tamhane2004superiority} targeted this problem with a test derived from
the union--intersection test of \citep{roy1953heuristic} and
the intersection--union test of \citep{berger1982multiparameter}.
We refer to \cite{wassmer1999procedures} for a more detailed
review of this area and \cite{lachin2014applications}
for recent advances. 

The hypothesis of interest in this paper, \eqref{H1},
is similar to but different from those considered in the above papers.
We investigate the direct application of the standard LRT to \eqref{H1} and
discover that a specific version of the LRT leads to a much improved
procedure that is particularly useful for our application.
We find a novel way to mildly relax
the size control to obtain a much more powerful test.
Simulation confirms that the new test retains tight control of
the type I error and is markedly more powerful than the LRT
as well as many of its competitors based on normal data.
The new method performs well in the context of 
monitoring multiple quality indices. 

The paper is organized as follows. In Section 2, we revisit
some basics of the LRT, introduce the new test, 
and review existing methods for normal data and one-sided
multi-parameter hypotheses. 
In Section 3, we give a brief background on the monitoring test
for forestry products and the application of the proposed method.
In Section 4, we present simulation results. We conclude in Section 5.

\section{Proposed and related methods}
The new approach was developed as a result of our observation
of the LRT under the normal model.
For this reason, we first quickly revisit the standard
likelihood approach and then introduce our approach.

\subsection{ LRT statistic}
Suppose we have an independent and identically distributed (iid) sample
$\bY_1,\ldots,\bY_n$ from a $p$-dimensional multi-normal distribution  
$\mbox{MVN}(\bmu,\bSigma)$. We first
consider the test problem for
\be
\label{H2}
H_0: \bmu \leq 0 \mbox{~~against~~} H_1:  \bmu \not \leq 0.
\ee

Let $\bX$ denote the sample mean $\bar{\bY}$ and 
\[
\bS
=
\frac{1}{n} \sum_{i=1}^n(\bY_i - \bar{\bY})(\bY_i - \bar{\bY})^T,
\]
a slightly altered sample variance.
It is well known that $\bX$ and $\bS$ together are complete and sufficient
for $\bmu$ and $\bSigma$ under the normal model. Hence, we may develop a likelihood-based method
as if they are the only observations.

After some simple algebra, the log-likelihood function is found to be
\[
\ell_n(\bmu,\bSigma)
=
-\frac{n}{2}\log\det(\bSigma) - \frac{n}{2}\tr\{\bSigma^{-1}[\bS+(\bX-\bmu)^T(\bX-\bmu) ]   \}.
\]
To develop an LRT, we search for the maximum point of $\ell_n(\bmu,\bSigma)$
under the null hypothesis and under the full model.
The solution under the full model is well known, with the unconstrained maximum
likelihood estimators of $\bmu$ and $\bSigma$ given by
\[
\hat{\bmu} = \bX; ~~ \hat{\bSigma} = \bS.
\]
This implies 
\[
\sup \ell_n(\bmu, \bSigma)  = - (n/2) \{ \log \det(\bS) + p\}.
\]

The solution under the null model is algebraically simple but slightly more abstract.
For each fixed $\bmu$, we find 
\[
\arg \max_{\bSigmas}  \ell_n(\bmu, \bSigma)
=  \hat{\bSigma}_{\bmu}
=
\bS+(\bX-\bmu)(\bX-\bmu)^T.
\]
This leads to the profile log-likelihood function of $\bmu$:
\ba
\ell_n(\bmu, \hat{\bSigma}_\mu)
&=&
-\frac{n}{2}\{\log\det(\bS+(\bX-\bmu)(\bX-\bmu)^T)+p   \}
\nonumber
\\
&=&
-\frac{n}{2}\{\log\det(\bS)+\log [1+(\bX-\bmu)^T\bS^{-1}(\bX-\bmu)]+p   \}.
\ea
The second equality is obtained by a linear algebra result
$\det(\mbox{\bf I} + \uu\vv^T )= 1+\uu^T\vv$ for any vector ${\uu}$ and $\vv$,
and by
\[
\bS+(\bX-\bmu)(\bX-\bmu)^T
= 
\bS^{1/2}\{ \mbox{\bf I}  +  [\bS^{-1/2}(\bX-\bmu)][ \bS^{-1/2}(\bX-\bmu)]^T\} \bS^{1/2} .
\]
Clearly, the profile likelihood is maximized if and only if
$(\bX-\bmu)^T\bS^{-1}(\bX-\bmu)$ is minimized with respect to $\bmu$
in the space of the null hypothesis. Let the solution to the minimization problem
be $\hat{\bmu}_0$. Geometrically, it is the projection of $\bX$ onto the null
space in terms of the Mahalanobis distance defined through the covariance matrix $\bS$.
Subsequently, we find the generic  expression of the LRT statistic:
\[
R_n 
=
2 \{ \sup \ell_n(\bmu, \bSigma) - \sup_{H_0} \ell_n(\bmu, \bSigma) \}
=n \log  \{ 1 + (\bX - \hat{\bmu}_0)^T \bS^{-1}  (\bX - \hat{\bmu}_0)\}.
\]
Note that $R_n$ is monotonic in 
\be
\label{ourT}
T_n = 
n (\bX - \hat{\bmu}_0)^T \bS^{-1}  (\bX - \hat{\bmu}_0).
\ee
Thus, the rejection region of the LRT statistic has the generic form
\be
\label{rejection}
C = \{ (\bY_1, \ldots, \bY_n): T_n  > c\}
\ee
for some $c$, which is called the critical value of the test.

By classical theory in mathematical statistics, if the size of the test is set to $\alpha$,
then the critical value $c$ will be chosen so that
\be
\label{LRT}
\sup_{\bmus \in H_0, \bSigmas>0} 
\pr\{ T_n > c; \bmu, \bSigma\} = \alpha
\ee
where we use $\pr(\cdot; \bmu, \bSigma)$ to indicate that the calculation is under
the $\mbox{MVN}(\bmu, \bSigma)$ distribution. 
According to \cite{perlman1969one}, the supremum is attained asymptotically
when $\bmu \to 0$ and $\bSigma$ approaches some singular matrix.
Specifically, he proved that for $H_0$ defined by \eqref{H1},
\ba
\sup_{\bmus \in H_0, \bSigmas>0} 
\pr\{  T_n > c; \bmu, \bSigma\} 
&=&
\frac{1}{2}\pr \big [\bF_{p-1, n-p+1} \geq \big ( \frac{1}{p-1} - \frac{1}{n} \big )c \big ]
\nonumber \\
&&
+
\frac{1}{2} \pr\big [ \bF_{p, n-p} \geq \big ( \frac{1}{p} - \frac{1}{n} \big )c \big ]
\label{eqn.LRTcriti}
\ea
where $\bF_{p, n}$ denotes an F-distributed random variable
with $p$ and $n$ degrees of freedom.
In other words, an LRT of size $\alpha$ will choose $c$ such that
\be
\label{LRTcrit}
\pr \big [\bF_{p-1, n-p+1} \geq \big ( \frac{1}{p-1} - \frac{1}{n} \big )c \big ]
+
\pr\big [ \bF_{p, n-p} \geq \big ( \frac{1}{p} - \frac{1}{n} \big )c \big ]
= 2 \alpha.
\ee

\subsection{Proposed test} 
The choice of $c$ in the LRT in \eqref{LRTcrit} ensures that
the type I error is at most $\alpha$ at any $(\bmu, \bSigma) \in H_0$.
When the dimension of the data $p=2$, the type I error is maximized
when $\bmu = \bf{0}$ and $\rho \to -1$ where $\rho$ is the
correlation coefficient. If the observations are from a distribution
with $\bmu = {\bf 0}$ and $\rho  = 0$, the type I error is far lower
than $\alpha$. In many applications, the user may be confident that
$\rho \geq 0$. If so, this choice is far too conservative. The size of
the test over the region of interest is much lower than the designated
$\alpha$. As a consequence, the power of the test is also much lower.

This consideration begs a question on the type I error of the test
at $\bmu = 0$ and a given $\bSigma$. 
Interestingly, an answer is readily available from \cite{nuesch1966problem}.
To state this result, we first introduce some notation.
When $\bX$ is $\mbox{MVN}(\bmu, \bSigma)$, we use the simplified notation
\[
\pr\{\bSigma\} = \pr\{ \bX> 0; \bmu = 0, \bSigma\}.
\]
Let $\cS$ be the collection of all nonempty subsets of $\{1, 2, \ldots, p\}$.
We use $\bX[i]$ for the $i$th entry of vector $\bX$.
For any $s \in \cS$, we use $\bX[s]$ for the subvector of $\bX$
consisting of components of $\bX[i]$ such that $i \in s$. 
Let $s'$ be the complement of $s$. 
With these, we use $\bSigma_s$ for the covariance matrix of $\bX[s]$ and
$\bSigma_{s'|s}$ for the covariance matrix of $\bX[s']$ conditional on $\bX[s]=0$.
We use the convention that when $s'$ is empty $\pr\{\bSigma_{s'|s}\} = 1$.
We use $|s|$ for the size of $s$.
In the following theorem, $T_n$ is the LRT statistic defined earlier. 

\begin{theorem}
\label{thm}
In the current setting, for any $c > 0$,
\[
\pr(T_n >c)
=
\sum_{s \in \cS}
\pr\big \{ \bF_{|s|, n - |s|} > \big ( \frac{1}{|s|} - \frac{1}{n}\big ) c \big \}
\pr\{\bSigma_s^{-1}\} 
\pr\{\bSigma_{s'|s}\}.
\]
\end{theorem}

In other words, the distribution of $T_n$ is a finite mixture of $F$-distributions.
The proof of this theorem is technically involved;
we refer to  \cite{nuesch1966problem} for the details.

The probabilities in the above theorem have generic analytical expressions
that can be found in \cite{kendall1941proof}. 
We are particularly interested in the case $p=2$. When $p=2$, without loss of
generality, we assume that $X$ has marginal variances $1$ and denote the
correlation coefficient as $\rho$. For $s$ such that $|s| = 1$, it is easy to
see that
\[
\pr\{\bSigma_s^{-1}\} = \pr\{\bSigma_{s'|s}\}= \frac{1}{2}.
\]
When $|s| = 2$, the correlationship coefficient specified by $\bSigma^{-1}$
is $- \rho$. Let $Z_1, Z_2$ be two independent $N(0, 1)$ random variables.
Then, $X_1 = Z_1$ and $X_2 = \sin (\gamma) Z_2 - \cos(\gamma) Z_1$ have
correlation $-\rho$ when $\gamma = \arccos(\rho)$ in the range of 0 and
$\pi$. Hence,
\[
\pr\{\bSigma_s^{-1}\} = \pr(Z_1 > 0; ~\sin (\gamma) Z_2 - \cos(\gamma) Z_1> 0)
=
\frac{\gamma}{2\pi}.
\]
In other words, we have
\be
\label{eq1.1}
\pr(T_n >c)
=
\frac{1}{2} \pr\big \{  \bF_{1, n-1}\geq \big (1 - \frac{1}{n} \big ) c\big \}
+
\frac{\arccos(\rho)}{2\pi} \pr\big \{  \bF_{2, n-2}\geq  \big (\frac{1}{2} - \frac{1}{n} \big )  c \big \}.
\ee
Consequently, if the value of $\rho$ is known and the observed value of $T_n$
is $t_{obs}$, we would have evaluated the $p$ value of the test to be
\[
\frac{1}{2} \pr \big \{  \bF_{1, n-1} \geq  \big (1 - \frac{1}{n} \big ) t_{obs}\big \}
+
\frac{\arccos(\rho)}{2\pi} \pr\big \{  \bF_{2, n-2} \geq  \big (\frac{1}{2} - \frac{1}{n} \big ) t_{obs}\big \} .
\]
This would lead to a much more powerful test than the classical
LRT. For instance, we would reject $H_0$ when $t_{obs} = 4.59$ when
$\rho$ is known to be $0$, while the LRT does not reject in this case. 
See Table \ref{tab_critical} for the critical values.
The LRT uses the critical value at $\rho = -1$, corresponding to the
least favorable distribution.

\begin{table}[ht]
\centering
\caption{Critical values of the LRT test when $\rho$ is known and $n=50, p=2$.}
\label{tab_critical}
\begin{tabular}{c|cccccc}
$\rho$ & $-1.0$ & $-0.9$ &$-0.5$ &  $0$ & 0.5   & $0.9$\\ \hline 
$c$     &  5.64   &   5.37   & 4.98   & 4.58 & 4.12  &  3.47\\
\end{tabular}
\end{table}


Motivated by the above discussion and calculations, we propose
a new test for $p=2$. First, we obtain the value of
$T_n$ and the sample correlation coefficient $\hat \rho$. 
With the observed value $t_{obs}$, we compute
\be
\label{pvalue}
\hat p
=
\frac{1}{2} \pr \big \{  \bF_{1, n-1} \geq  \big (1 - \frac{1}{n} \big ) t_{obs}\big \}
+
\frac{\arccos(\hat \rho)}{2\pi} \pr\big \{  \bF_{2, n-2} \geq  \big (\frac{1}{2} - \frac{1}{n} \big ) t_{obs}\big \} .
\ee
The test rejects $H_0$ when $\hat p < \alpha$, where
$\alpha$ is the designated size of the test.

Our idea is not limited to $p = 2$. The 
analytical form of $\hat p$ (the p-value of the test)
is more complex in the general case
but can be calculated according to Theorem \ref{thm}.
We do not present the details here since the interested user can
work them out with some algebraic effort. We call the
new test the mLR test.

The type I error of the mLR test may in theory exceed $\alpha$ at some
specific $\rho$ values very close to $-1$. Our simulation
experiments show that the degree of inflation is negligible.

\subsection{Application to non-normal data}
In applications, the data are often collected from non-normal populations.
Nevertheless, it is generally possible to obtain a good estimate of the vector
parameter $\btheta$ of dimension $p$ and its covariance matrix.
We consider the situation where
\[
\sqrt{n} \bS_n^{-1/2} (\hat{\btheta} - \btheta) \to \mbox{MVN}(0, \bf{I})
\]
in distribution when some index, likely the sample size $n$, goes to
infinity.

Suppose it is of interest to test the hypothesis in the form of \eqref{H1}
and, without loss of generality, $\btheta^* = {\bf 0}$.
The proposed modified LRT can be applied to this problem
by setting $\bX = \hat{\btheta}$ and $\bS = \bS_n$.
The computation of $T_n$ and $\hat p$ can then be carried out
in the same way. We reject the null hypothesis when $\hat p < \alpha$.
When the sample size $n$ is large, one may use $\chi_p^2$ to replace
$\bF_{p, n}$ and so on to give an approximate $\hat p$.

\subsection{Other methods}
As pointed out earlier, there exist many methods to handle the hypothesis
test problem under a multivariate normal model. It is helpful to see how
the proposed method differs. For brevity, we give a quick introduction
to just two methods. We still assume that an iid sample $\bY_1, \ldots, \bY_n$
from $\mbox{MVN}(\bmu, \bSigma)$ is given and will continue to use
some of the notation introduced earlier.

\vspace{1ex}
\noindent
{\bf {Union--Intersection Test}}
In the union--intersection test (UIT),  
we start by defining  subnull hypotheses
$H_{0j} = \{ \bmu: ~\mu_j \leq 0 \}$ for $j=1, 2, \ldots, p$. 
Clearly, $H_0=\bigcap_{j=1}^p H_{0,j}$.
This means that if any $H_{0,j}$ is false, then $H_0$ is also false. 
Thus, one may test the validity of $H_{0,j}$ for each $j$. 
We reject $H_0$ if any $H_{0,j}$ is rejected.

When $\bSigma$ is known to be ${\bf I}$, we may reject $H_{0,j}$ when the component sample mean
of the $j$th component $\bar{\bY}_j >c$ for some critical value $c>0$. 
We reject $H_0$ when 
\[
\max\{\bar{\bY}_j: ~ j = 1, \cdots, p \}>c.
\]
Note that under the null hypothesis
\ba
\pr(\max\{\bar{\bY}_j, ~j =1, \cdots, p\} > c )
&=&
1-\pr(\max\{\bar{\bY}_j, ~j =1, \cdots, p\} < c)
\nonumber \\
&=&
1-\prod_{j=1}^p \pr(\bar{\bY}_j < c). 
\ea
Hence, we may choose  $c = z_{(1-\alpha)^{1/p}}/\sqrt{n}$ to obtain
a size $\alpha$ test, where $z_{(1-\alpha)^{1/p}}$ is
the lower $(1-\alpha)^{1/p}$ quantile of the standard normal distribution.

When $\bSigma$ is unknown, we may conduct a one-sided
$t$-test of size $\alpha/p$ for $H_{0j}$ for $j=1, 2, \ldots, p$.
We reject $H_0$ when any $H_{0j}$ is rejected. 
By the Bonferroni inequality we see that the size of this test below $\alpha$.
It is well known that a test formed by Bonferroni correction
tends to be very conservative.

\vspace{1ex}
\noindent
{\bf PW test.}
\cite{perlman2003validity} were among the first to take note of the
conservative nature of both UIT and LRT. 
In particular, they suggested that the boundary of $H_{0}$ 
can be decomposed into subspaces of varying dimensions.
For instance, when $p=2$,  the boundary of $\{\bmu \leq 0\}$ 
is decomposed into
\[
B_1= \{\mu_1 = 0, \mu_2 = 0\},~
B_2= \{\mu_1 < 0, \mu_2 = 0\},~
B_3= \{\mu_1 = 0, \mu_2 < 0\}.
\]
The dimension of $B_1$ is $0$ and that of $B_2$ and $B_3$ is $1$.
If the sample mean $\bX \in H_0$, then $T_n = 0$.
Otherwise, the maximum of the distances from $\bX$ to $B_1$, $B_2$, or $B_3$
is taken as $T_n$. 
The information on the source of the maximum is then discarded,
and the size of $T_n$ is measured against the least favorable distribution,
which corresponds to $\bmu \in B_1$ and $\rho = -1$. 

Perlman \& Wu fix the conservative nature of the LRT by having different
critical values depending on the location of $\bX$ with respect to
$B_1$, $B_2$, or $B_3$. 
Let
\begin{eqnarray*}
M_1 &=& \{n\bX^T \bS^{-1} \bX > c_{2, \alpha}\}\\
M_2 &=& \{\frac{\bX_1}{\sqrt{\bS_{11}/(n-1)}} > t_{n-1, \alpha}\}\\
M_3 &=& \{\frac{\bX_2}{\sqrt{\bS_{22}/(n-1)}} > t_{n-1, \alpha}\},
\end{eqnarray*}
where $c_{2,\alpha}$ is the critical value of the LRT test
of size $\alpha$, according to \eqref{eqn.LRTcriti}, and $\bS_{11}$ and
$\bS_{22}$ are entries of matrix $\bS$.
The PW test \citep{perlman2006some} rejects $H_0$ 
when $\bX \in M_1 \cap (M_2 \cup M_3)$.
That is, $H_0$ is rejected when $B_1$ is rejected and one of 
$B_2$ and $B_3$ is also rejected.

\begin{figure}[ht]
\centering	
	\includegraphics[width=0.6\textwidth]{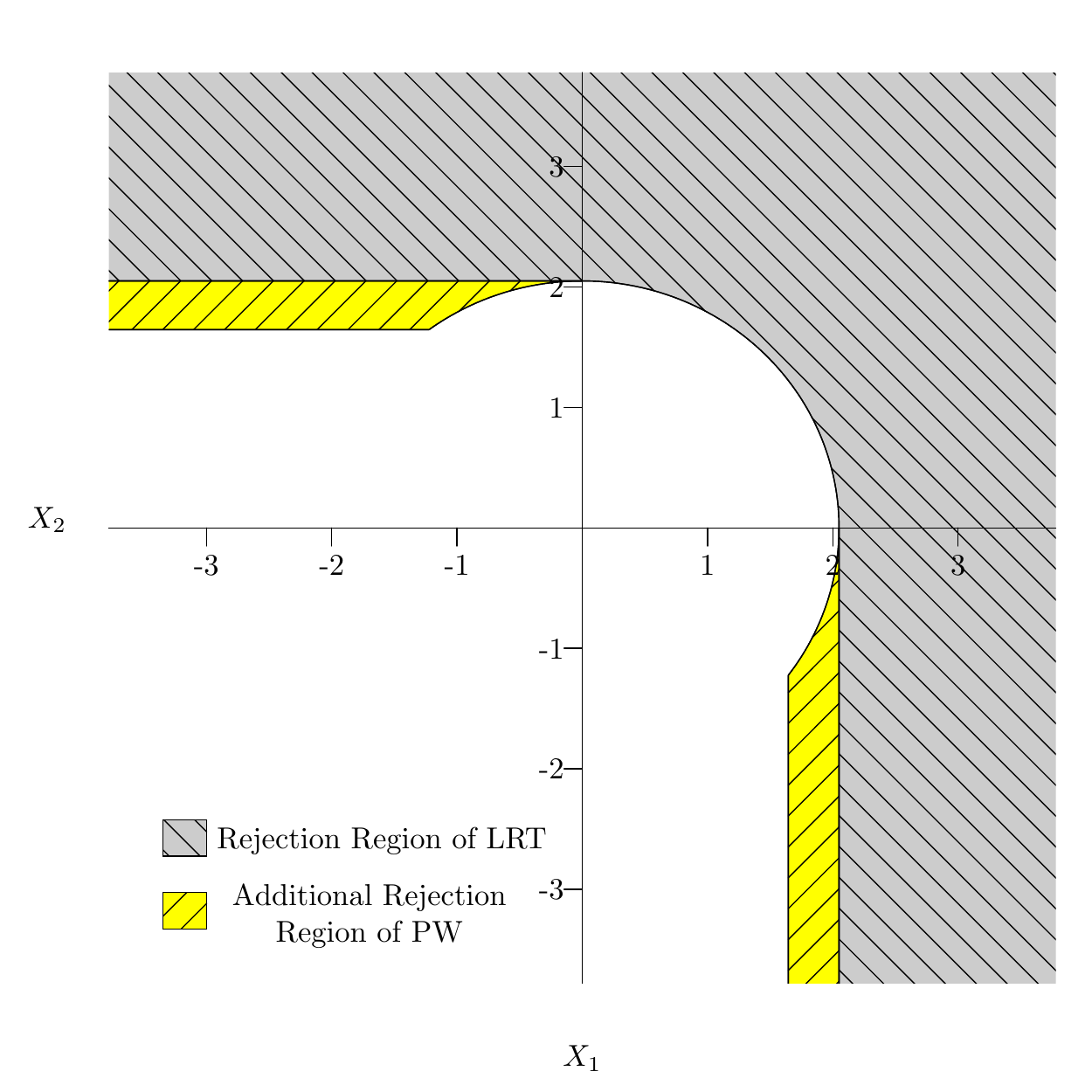}
	\caption{Rejection regions of LRT and PW test.}
\label{fig_bound3}
\end{figure}

We can verify that the rejection region of the PW test
covers the rejection region of the LRT; see Figure \ref{fig_bound3}. 
At the least favorable distribution where  $\rho = -1$, its type I
error will exceed $\alpha$, as is the case for our method.
When $\rho=-0.9$ the type I error of the PW test is $5.46\%$ based on
our simulations.


\section{Application to monitoring test}
 \label{monitor}

The proposed modified LRT is developed with an application in mind.
As discussed by \cite{verrill2015simulation}, forestry is concerned with
monitoring the lower quantiles of the mechanical strength distribution.
Many researchers focus on the $5$th quantile. In this paper, we simultaneously monitor
several quality parameters of the mechanical strength distribution.
In this section we demonstrate the usefulness of the modified LRT.

The modified LRT may be used  in many ways and many applications.
We, however, focus on the specific setting and inference
methods developed in \cite{chen2016monitoring}. We refer to this paper for 
more detailed background information but provide some necessary
description of the data and inference methods here.

The data under consideration are assumed to be a random sample
from $m+1$ populations with some clustered structure:
\[
\yy_{k, j}^T = (y_{k,j,1}, \ldots, y_{k,j,d}): k = 0, 1, \ldots, m;~ j = 1, 2, \ldots, n_k.
\]
In this setting, $k$ is the identity of the population, $d$ is the cluster size,
and $n_k$ is the number of clusters sampled from the $k$th population.

Let $F_k(\yy)$ be the cumulative joint distribution (CDF) of $\yy_{k,j}$. 
The nature of the data implies that $F_k$ is exchangeable. 
The exchangeability implies an identical marginal distribution,
which will be denoted $G_k(y)$. The target of the monitoring test
is hence $G_k(y)$. We wish to be alerted when
$G_k(y)$ is stochastically smaller than $G_0(y)$ in some respect.
As pointed out earlier, we may test if
$G_k$ is lower than $G_0$ in the 5\% quantile or
the median.

Because the $G_k$'s are of a similar nature,
\cite{chen2016monitoring} suggested that the density ratio model (DRM) 
\citep{anderson1979multivariate} is appropriate.
Specifically, they assumed that these distributions
are related through the following equation:
\begin{equation}
\label{eqn_DRM}
	\frac{\diff G_k(y)}{\diff G_0(y)}=\exp\{\bbeta^T_k \bq(y)\}
\end{equation}
for a suitably selected function $\bq(y)$ of dimension $q$ 
with unknown parameter vectors $\bbeta_k$.

Based on the DRM, \cite{chen2016monitoring} 
proposed the following composite empirical likelihood (EL):
\begin{equation}
\label{eqn_CEL}
	L(G_0,\bbeta)
	=\prod_{k,j,l}\{\diff G_k(y_{k,j,l})\}
	=\Big(\prod_{k,j,l} p_{k,j,l}\Big) 
	\exp \Big\{\sum_{k,j,l}\bbeta^T_k \bq(y_{k,j,l})\Big\}
\end{equation}
where  $G_0(y)=\sum_{k,j,l} p_{k,j,l}1(y_{k,j,l}\leq y)$.
The DRM assumption implies
\begin{equation}
\label{eqn_condition}
	\sum_{k,j,l} p_{k,j,l} \exp\{\bbeta^T_r\bq(y_{k,j,l})\} = 0
\end{equation}
for $r=0, 1, \ldots, m$.

Some algebra shows that the above composite EL
has a dual form:
\begin{equation}
\label{eqn_DEL}
	\ell_n(\bbeta)
	=
	- \sum_{k,j,l} \log [\sum_{r=0}^m \rho_r \exp\{\bbeta_r^T \bq(y_{k,j,l})\}]
	+ \sum_{k,j,l} \bbeta_k^T \bq(y_{k,j,l}).
\end{equation}
Many of the numerical computations are done via the dual form. 

Let the maximum composite EL estimator be 
$\hat{\bbeta} = \arg \max_\beta \ell_n(\bbeta)$.
Let
\[
\hat{G}_r(y)
=
\sum_{k,j,l}  \hat{p}_{k,j,l} 
\exp \{ \bbetahat_r^T \bq(y_{k,j,l}) \} \ind ( y_{k,j,l} < y)
\]
be the fitted CDF, with the obvious notation $\hat{p}_{k,j,l}$.
By the invariance property of the maximum likelihood estimation,
we estimate the population means and quantiles by
\[
\muhat_r
=
\sum_{k,j,l} \hat{p}_{k,j,l}\exp \{\bbetahat_r^T \bq(y_{k,j,l})\} y_{k,j,l}
\]
and
\[
\hat{\xi}_r
=
\hat{\xi}_{r, \alpha} = \inf \{y: \hat{G}_r(y) \geq \alpha\}
\]
where $\alpha$ denotes the level of the quantile.
It has been shown that the parameter estimators
are asymptotically normal. For instance, in obvious notation,
\[
\sqrt{n} \{(\hat \xi_1, \hat \xi_2) -  (\xi_1, \xi_2)\} \to N(0, \bSigma).
\]
A cluster-based bootstrap method proposed by
\cite{chen2016monitoring} can be used for the consistent estimation of $\bSigma$.

We are now ready to apply the modified LR test to the
one-sided test problem for multiple parameters.
Suppose $\btheta$ is a vector-valued parameter. Let $\hat{\btheta}$
be its MLE and $\bS^*$ be its bootstrap variance estimator 
given in \cite{chen2016monitoring}. The monitoring test problem is transformed to
the problem of testing for some hypothesis in the form of
\eqref{H1}. 
When 
\[
\btheta = (\xi_{1, 0.05} - \xi_{0, 0.05}, \xi_{1, 0.50} - \xi_{0, 0.50})^T,
\]
testing for \eqref{H1} involves monitoring whether $G_1$ has 
simultaneously maintained the $5$th percentile and the median
of  the wood strength distribution compared to $G_0$.
In the presence of multiple populations, the test is more efficient
if we also utilize information from $G_2$, $G_3$, and so on
\citep{chen2016monitoring}. Depending on the monitoring target, 
other forms of $\btheta$ can easily be specified.

The null hypothesis of interest is $\btheta \geq 0$.
To apply the proposed modified LRT, we compute the value
of $T_n$ given in \eqref{ourT} with 
\[
\bX = - \hat \btheta; ~~ \bS = n\bS^*.
\]
The reason for the negative sign in $\bX = - \hat \btheta$
is to reconcile the opposite inequalities specified
in \eqref{H1} and \eqref{H2}. 
We compute the p-value of the test according to \eqref{pvalue}.
Clearly, we could as easily use other tests based on $\bX$ and $\bS$.

\section{Simulation and example}

In this section, we use simulation to discover the pros and cons
of three tests: LRT, PW, and the proposed mLR
for one-sided hypotheses. 
We do not include UIT because this method has been shown
to be inferior by \cite{perlman2003validity} 
and \cite{perlman2006some}.
As pointed out earlier,
the type I errors of the mLR and PW tests likely exceed the
desired size for some distributions. It is important to explore
how serious the errors become and the features of the corresponding distributions.

We focus on the situation where the dimension of the parameter
$p=2$ with a sample of size $n=50$ from various multivariate
normal distributions. 

\subsection{Multivariate normal samples}
It can easily be seen that the test problem of interest is invariant to
the variance of the marginal distributions. When $p=2$, this implies
that we need consider only the covariance matrices in the following
form:
\[
\bSigma
=
\begin{pmatrix}
1 & \rho\\ 
\rho & 1
\end{pmatrix}.
\]
We generated data from $4 \times 5$ null models with a range of
correlation coefficients:
\[
\rho = -0.9, -0.5, 0.0, 0.5, 0.9.
\]
From each model, we generated 100,000 samples of size $n=50$.
We set the nominal rejection rate, or size of the test, to $0.05$.
The values of the population mean $\bmu$ and the percentage of times 
when the null hypothesis is rejected by these four tests are summarized
in Table \ref{tab1}.

{\bf Null models.}
Let us first examine the results for $\bmu = (0, 0)^T$ at which the
null hypothesis is true.
The results in Table \ref{tab1} support the theory that LRT and UIT tightly 
control the type I error. However, they achieve this goal by being very
conservative at $\rho = 0.5, 0.9$.
The PW test improves on LRT and UIT in terms of being less conservative,
but at the cost of exceeding the nominal level at $\rho = -0.9$.
The type I errors of the proposed mLR over this range of $\rho$
are very close to the nominal level. 

When $\bmu$ goes from $(0, -0.1)$ to $(0, -0.3)$, the null hypothesis
remains true. Since it makes the model move toward the ``interior'' of $H_0$,
the type I errors of these tests become lower, as expected.

{\bf Alternative models.}
We also carry out simulation for three sets of alternative distributions.
In the first, both marginal means become greater than 0 at the same rate.
In the second, just one of the marginal means becomes greater than 0.
In the third, two marginal means move in opposite direction.
The simulated powers of the three tests are given in
the second, third and fourth blocks of Table \ref{tab1}.

Clearly, LRT has lower power than PW and mLR for the alternative distributions. 
The comparison between PW and mLR is not clear-cut: mLR is uniformly more powerful 
than PW for the first set of alternative distributions (second block of  Table \ref{tab1}).
For the second set (third block of Table \ref{tab1}) mLR has higher power than PW 
when  $\rho = -0.9, -0.5$, and $0$; comparable
power when $\rho = 0.5$; and slightly lower power when $\rho = 0.9$.
For the third set (fourth block of Table \ref{tab1}) PW is more powerful.

Based on the simulation results, we recommend using the PW test in applications
where the two quality indices may move in opposite directions.
If the two indices are likely to move in the same direction, 
mLR is preferable.


\begin{table}[p]
	\centering
	\caption{Type I errors for one-sided tests (\%)}
		\label{tab1}
	\setlength{\tabcolsep}{0.4em}
	\small{
	\begin{tabular}{|c|ccc|ccc|ccc|ccc|ccc|}
		\hline
		& \multicolumn{3}{c|}{$\rho=-0.9$} & \multicolumn{3}{c|}{$\rho=-0.5$} & \multicolumn{3}{c|}{$\rho=0$}  & \multicolumn{3}{c|}{$\rho=0.5$} & \multicolumn{3}{c|}{$\rho=0.9$} \\\hline
$\bmu$	  & LRT & PW & mLR & LRT & PW & mLR & LRT & PW & mLR & LRT & PW & mLR & LRT & PW & mLR \\ 
		\hline
(0,0) & 4.49 & 5.46 & 5.08 & 3.79 & 4.79 & 5.05 & 3.05 & 4.03 & 4.98 & 2.46 & 3.43 & 5.00 & 1.66 & 2.64 & 4.94 \\ 
(0,-.1) & 1.34 & 4.23 & 1.56 & 1.54 & 3.03 & 2.12 & 1.62 & 2.91 & 2.69 & 1.34 & 2.70 & 2.87 & 1.18 & 3.62 & 3.65 \\ 
(0,-.2) & 1.23 & 5.13 & 1.40 & 1.11 & 3.60 & 1.56 & 1.20 & 3.32 & 2.04 & 1.20 & 3.69 & 2.59 & 1.15 & 4.90 & 3.55 \\ 
(0,-.3) & 1.14 & 5.12 & 1.31 & 1.14 & 4.48 & 1.58 & 1.21 & 4.28 & 2.05 & 1.18 & 4.50 & 2.55 & 1.13 & 4.89 & 3.45 \\ 
\hline
(.1,.1) & 84.4 & 84.4 & 85.7 & 27.6 & 27.9 & 32.4 & 16.1 & 17.1 & 22.3 & 11.2 & 12.9 & 18.8 & 7.74 & 10.2 & 17.3 \\ 
(.2,.2) & 100 & 100 & 100 & 75.8 & 75.8 & 80.1 & 46.8 & 47.3 & 56.3 & 32.4 & 34.1 & 45.3 & 23.6 & 27.5 & 40.8 \\ 
(.3,.3) & 100 & 100 & 100 & 97.9 & 97.9 & 98.6 & 80.3 & 80.4 & 86.2 & 62.2 & 63.1 & 74.3 & 49.4 & 53.4 & 68.7 \\ 
(.4,.4) & 100 & 100 & 100 & 99.9 & 99.9 & 99.9 & 96.3 & 96.3 & 97.9 & 86.2 & 86.5 & 92.5 & 75.1 & 77.6 & 88.3 \\ 
		\hline
(0,.1) & 33.9 & 34.1 & 36.0 & 12.4 & 13.4 & 15.5 & 8.88 & 10.6 & 13.0 & 6.79 & 9.64 & 12.4 & 5.72 & 12.6 & 13.5 \\ 
(0,.2) & 86.1 & 86.1 & 87.3 & 33.4 & 34.3 & 38.8 & 23.1 & 26.3 & 30.4 & 19.4 & 27.5 & 29.6 & 18.4 & 38.9 & 33.9 \\ 
(0,.3) & 99.5 & 99.5 & 99.6 & 62.7 & 63.1 & 68.0 & 46.7 & 50.8 & 55.9 & 42.3 & 56.0 & 55.6 & 41.9 & 67.2 & 61.3 \\ 
(0,.4) & 100 & 100 & 100 & 86.4 & 86.6 & 89.3 & 72.0 & 75.3 & 79.2 & 68.8 & 82.0 & 79.7 & 68.2 & 87.3 & 83.6 \\ 
\hline
(-.1,.1) & 8.19 & 10.3 & 9.04 & 7.07 & 9.89 & 9.14 & 6.26 & 10.1 & 9.44 & 5.83 & 11.7 & 10.6 & 5.53 & 16.7 & 13.4 \\ 
(-.2,.2) & 21.1 & 25.9 & 22.7 & 19.6 & 27.8 & 23.6 & 18.7 & 31.4 & 25.2 & 18.4 & 37.5 & 28.4 & 18.6 & 40.5 & 34.4 \\ 
(-.3,.3) & 44.0 & 51.2 & 46.1 & 42.1 & 55.9 & 47.5 & 41.8 & 62.8 & 50.9 & 41.7 & 66.8 & 54.8 & 41.6 & 67.1 & 61.3 \\ 
(-.4,.4) & 69.3 & 76.2 & 71.2 & 68.5 & 81.8 & 73.2 & 68.3 & 86.2 & 75.9 & 68.3 & 87.4 & 79.3 & 68.5 & 87.4 & 83.7 \\ 
\hline
	\end{tabular}
}
\end{table}

\subsection{Application to multiple quality indices in monitoring context}
We now study the use of the proposed test for multi-dimensional quality
indices in monitoring. We simulate data with a cluster structure, as discussed in Section \ref{monitor}.
We compare the LRT and PW test and we again omit UIT.

We consider the situation where clustered random samples from 
$m+1 = 4$ populations are available and the cluster size $d = 5$.
We use $B = 999$ bootstrap repetitions for the variance estimation.
To paint a more complete picture, we simulated data from two clustered population sets:
one is multivariate normal and the other is multivariate gamma. The reliability literature 
indicates that these are sensible
models for data from quality indices. We emphasize
that the data analysis does not assume knowledge of the data-generating distributions.

\paragraph{Multivariate clustered normal populations}
We first perform simulation by generating individual response values
from the following random effect model:
\[
y_{k,j,l}
=
\mu_k + \gamma_{kj} + \epsilon_{kjl}.
\]
In the wood product application, $y_{k, j, l}$ is the mechanical strength of
a piece of wood from the $k$th population, $j$th cluster, and $l$th unit.
We generate $\gamma_{kj}$ from $N(0, \sigma^2_{\gamma, k})$. 
Since $\gamma_{kj}$  is shared by all the units in cluster $j$ in the $k$th population,
it induces within-cluster positive correlation. 
We generate $\epsilon_{kjl}$ from $N(0, \sigma^2_e)$, which reflects the noise
in the mechanical strength.
The marginal distributions $G_k$ are all normal, but this fact will
not be used in the hypothesis test. Instead, we use DRM with
$\bq(y) = (1, y, y^2)^T$.

The problem of interest in the targeted application is whether or not the
$5$th percentile and the median of the mechanical strength of year
$k > 0$ are maintained compared to some base year $k=0$.
Let $\xi_{k, \alpha}$ be the $\alpha$th percentile of $G_k$.
Let 
\[
\btheta_k = (\xi_{k, 0.05} - \xi_{0, 0.05};~ \xi_{k, 0.50} - \xi_{0, 0.50}).
\]
For the purposes of illustration, we test, for each $k=1, 2, 3$ not simultaneously,
\[
H_0: \btheta_k \geq 0
\mbox{ against the alternative }
H_a: \btheta_k \not \geq 0.
\]
Clearly, the proposed test can be used for any other suitable
quality indices. The same is true for the LRT and the PW test.

The simulation was conducted with three sets of parameters:
\[
\begin{array}{|c|c|c|c|c|}
\hline
 & (\mu_0, \cdots, \mu_3) & (\sigma_{\gamma, 0}, \ldots, \sigma_{\gamma,3}) & \sigma_e & \mbox{Feature}\\
 \hline
\mbox{I} &  (15.5, 15.5, 14.7, 14.0) & (1.2, 1.2, 1.0, 1.0)   & 2.0 & \mbox{$\xi_{.05}, \xi_{.50}$ reduced}   \\
\hline
\mbox{II}  &  (15.5, 15.2, 15.0, 14.7) & (2.0, 1.794, 1.653, 1.436) &1.0 & \mbox{$\xi_{.05}$ reduced}  \\
\hline
\mbox{III} & (15.5,15.5,15.5,15.5) & (1.0,1.2,1.4,1.6) & 1.0 & \mbox{$\xi_{.50}$ reduced}\\
\hline
\end{array}
\]
The numbers of clusters are chosen to be $(n_0, n_1, n_2, n_3) = (25, 30, 40, 40)$.
The quantile and median values are given by
\[
\begin{array}{|c|c|c|}
\hline
& (\xi_{0,0.5}, \cdots, \xi_{3,0.5}) &  (\xi_{0,0.05}, \cdots, \xi_{3,0.05}) \\
\hline
\mbox{I} &  (15.50, 15.50, 14.70, 14.00) & (11.66, 11.66, 11.02, 10.32)   \\
\hline
\mbox{II}  &  (15.50, 15.20, 15.00, 14.70) & (11.82, 11.82, 11.82, 11.82) \\
\hline
\mbox{III} & (15.50,15.50,15.50,15.50) & (13.17, 12.93, 12.67, 12.40)  \\
\hline
\end{array}
\]


In the first setting, the first two populations are identical and the other two populations
have a lower $5$th percentile and median.
This arrangement allows us to investigate the type I error by testing
$\btheta_1 \geq 0$ and the power for $\btheta_2 \geq 0$ and $\btheta_3 \geq 0$.
In the second setting, the four populations have the same median, but
the $5$th percentile reduces from the first to the last population.
In the third setting, the four populations have the same $5$th percentile, but
the median reduces from the first to the last population.

We set the number of repetitions to $10,000$.
The simulated rejection rates for the three hypotheses are 
summarized in Table \ref{tab_cluster_test1}.

\begin{table}[ht]
	\centering
		\caption{Simulated rejection rates for normal data (\%)}
		\label{tab_cluster_test1}
	\begin{tabular}{|c|rrr|rrr|rrr|}
		\hline
		&\multicolumn{3}{c|}{\mbox{Setting I}} &\multicolumn{3}{c|}{\mbox{Setting II}}& \multicolumn{3}{c|}{\mbox{Setting III}}\\
				\hline
		$H_0$ & LRT & PW  & mLR & LRT & PW & mLR & LRT & PW & mLR \\ 
		$\btheta_1\geq 0$ & 2.93 & 3.86 &  5.91 & 4.50 & 5.20  & 8.20 & 6.29 & 9.34  & 12.08\\
		$\btheta_2\geq 0$ & 47.35 & 52.55  & 62.01& 7.20 & 10.30  & 13.30 & 16.45 & 25.55 & 25.87 \\ 
		$\btheta_3\geq 0$ & 95.81 & 96.83  & 98.20 & 14.30 & 22.20 & 24.70 & 29.18 & 44.44 & 42.93\\  
		\hline
	\end{tabular}
\end{table}

Recall that in Setting I, the null hypothesis $\btheta_1 \geq 0$ is true. The simulation results clearly
show that the faithful LRT has a much lower type I error than the nominal size of 5\%. This is not bad in
itself. The problem is that the lower type I error is at the cost of a much lower power for rejecting
$\btheta_2 \geq 0$ and $\btheta_3 \geq 0$ compared to the other methods.
Comparing PW and mLR shows that PW is also too conservative and therefore
has low power. The mLR has higher power but also higher type I error. 

The null hypotheses for Settings II and III are false, and so power is measured by the rejection 
of the hypothesis.
The simulation results in Table \ref{tab_cluster_test1} generally favor mLR. Overall, we
conclude that the proposed mLR works well.

\paragraph{Multivariate clustered gamma populations}
We now perform simulation by generating individual response values
from multivariate clustered gamma populations.

One way to create multivariate clustered gamma observations is as follows.
Let $U_1,\ldots,U_d$ be $d$ iid random variables following beta distributions 
with shape parameters $a$ and $b$. 
Further, let $W$ be a gamma-distributed random variable with shape parameter $a+b$ 
and rate parameter $\beta$. 
Then 
\[
\bY=W (U_1,\ldots,U_d)^T
\]
is multivariate gamma $MG(a,b,\beta)$ 
with correlation $\mbox{cor}(Y_i, Y_j)=a/(a+b)$ for all $1 \leq i < j \leq d$. 
The marginal distribution of $\bY_1=U_1W$ is gamma with shape parameter $a$ and 
rate parameter $\beta$. When $b=\infty$, $Y_1,\ldots,Y_d$ become independent;
see \cite{nadarajah2006some}.

The simulation was conducted with three sets of parameters:
\[
\begin{array}{|c|c|c|c|c|}
\hline
 & (a_0, \cdots, a_3) & (\beta_0, \ldots, \beta_3) & b & \mbox{Feature}\\
 \hline
\mbox{I} &  (8.0, 8.0, 7.0, 6.0) & (1.00,1.00,1.05,1.10)   & 14 & \mbox{$\xi_{.05}, \xi_{.50}$ reduced}   \\
\hline
\mbox{II}  & (8.0, 8.5, 9.0, 10) & (1.00,1.09,1.18,1.36) &14 & \mbox{$\xi_{.05}$ reduced}  \\
\hline
\mbox{III} & (8.0,7.0,6.0,5.0) & (1,0.87,0.74,0.61) & 14& \mbox{$\xi_{.50}$ reduced}\\
\hline
\end{array}
\]
The quantile and median values are given by
\[
\begin{array}{|c|c|c|}
\hline
& (\xi_{0,0.5}, \cdots, \xi_{3,0.5}) &  (\xi_{0,0.05}, \cdots, \xi_{3,0.05}) \\
\hline
\mbox{I} &  (7.67, 7.67, 6.35, 5.15) & ( 3.98, 3.98, 3.13, 2.38)   \\
\hline
\mbox{II}  &  (7.67, 7.49, 7.35, 7.11) & (3.98, 3.98, 3.98, 3.98) \\
\hline
\mbox{III} & (7.67, 7.67, 7.67, 7.67) & (3.98, 3.78, 3.53, 3.23)  \\
\hline
\end{array}
\]

We test the same hypotheses as for the multivariate clustered normal populations.
The results are given in Table  \ref{tab_cluster_test2}.
\begin{table}[ht]
	\centering
		\caption{Simulated rejection rates for gamma data (\%)}\label{tab_cluster_test2}
	\begin{tabular}{|c|rrr|rrr|rrr|}
		\hline
		&\multicolumn{3}{c|}{\mbox{Setting I}} &\multicolumn{3}{c|}{\mbox{Setting II}}& \multicolumn{3}{c|}{\mbox{Setting III}}\\
				\hline
		$H_0$ & LRT & PW  & mLR & LRT & PW & mLR & LRT & PW & mLR \\ 
		$\btheta_1\geq 0$ & 2.79 & 3.76 & 5.69 & 76.07 & 77.48 & 86.23 & 99.99 & 99.99 & 100.0 \\
		$\btheta_2\geq 0$ & 4.17 & 5.43 & 7.96 & 6.25 & 8.96 & 12.30& 13.80 & 21.82 & 23.51 \\ 
		$\btheta_3\geq 0$ & 6.01 & 8.72 & 11.27 & 14.21 & 21.17 & 22.84 & 32.87 & 47.61& 45.44\\  
		\hline
	\end{tabular}
\end{table}

Our observations are similar to those for the multivariate clustered normal
populations. Both LRT and PW are too conservative: the type I error is much lower than 5\%
in Setting I, for the null hypothesis $\btheta_1 \geq 0$. 
The PW test is also too conservative and therefore has low power. 
The mLR has higher power but also higher type I error. 
The overall impression is that the proposed mLR works well.

\subsection{Data analysis}
We now apply our method to a real forestry data set. 
It contains 398 modulus of rupture (MOR) measurements from In-Grade samples 
and 408 MOR measurements from monitoring samples obtained in 2011/2012.
Both \cite{chen2016monitoring} and \cite{verrill2015simulation} found that the 5th quantile
is markedly reduced in the monitoring sample with high statistical significance.
\begin{table}[ht]
	\centering
	\caption{Sample quantiles of forestry data}
	\begin{tabular}{rrr}
		\hline
		& 5\% & 50\% \\ 
		\hline
		In-Grade & 2.64 & 5.28 \\ 
		2011/2012 & 1.87 & 3.71 \\ 
		\hline
	\end{tabular}
\end{table}
We certainly expect that any one-sided hypothesis tests for the 5th quantile and
the median of MOR will produce a statistically significant outcome.
In this analysis, we used the basis function $\bq(y) = (1, y, y^2, \log y)$
suggested by \cite{chen2016monitoring}.
The estimated differences in the 5th quantile and
the median are $(\hat \theta_{0,1;0.05}, \hat \theta_{0,1;0.5})= (- 0.69, - 1.53)$. 
By the bootstrap method recommended by \cite{chen2016monitoring},
the asymptotic covariance matrix of this estimator is estimated as
\[
\bS_n
=
\begin{pmatrix} 
0.01282 & 0.01586\\ 
0.01586 &  0.04022
\end{pmatrix}.
\]
We now use $\bX = ( 0.69, 1.53)^T$ and $\bS = n \bS_n$
to compute $T_n$ defined in \eqref{ourT}.
We find $T_n = 59.3$ and $\hat p =2.30 \times 10^{-14}$ by \eqref{pvalue}.
Hence, the null hypothesis is rejected with strong
statistical evidence.

Note that the estimated correlation coefficient is $\hat \rho = 0.70$
in this example. This is the value used to compute $\hat p$.
When the LRT is applied to this problem, we compute
the p-value as if $\rho = -1$, giving
$7.15\times 10^{-14}$. The result remains sufficiently significant,
but there is a large drop in the level of significance.
The p-value of the PW test is the same in this case.

The two populations in this example are so different that
the quality deterioration is detected by any reasonable methods. 
To demonstrate more subtle differences between methods, 
we artificially inflate every data point of the 2011/2012 sample by a factor of 1.35. 
The two samples now have much closer sample-quality indices:
the estimated differences in the 5th quantile and
the median are $(\hat \theta_{0,1;0.05}, \hat \theta_{0,1;0.5}) = (-0.166, -0.009)$. 
The estimated asymptotic covariance matrix of this estimator is
\[
\bS_n
=
\begin{pmatrix} 
0.0081& 0.0156\\
0.0156& 0.0545
\end{pmatrix}.
\]
We now find $T_n= 3.41$, and
the p-values based on LRT, PW, and mLR are 
$0.123$, $0.032$, and $0.053$. 
Because the change in the median is so small, the PW test arrives
at its p-value primarily because of the large $|\hat \theta_{0,1; 0.05}|$.
In comparison, mLR takes a more balanced view of the two indices, and
the differences in the median and 5\% quantile between the
two populations are judged not significant at the 5\% level. The LRT is too conservative,
as our simulations predicted.

\section{Conclusions}
One-sided multi-parameter hypothesis tests arise in many applications, and
there are many effective test methods under normal models with a solid theoretical basis. 
We are particularly interested in testing 
whether two quality indices are reduced over time. 
The existing methods have room for further
improvement, particularly in the context of our application.
We propose a new test for this context.
In particular, we have developed a strategy for applying the method to general
one-sided multi-parameter hypotheses.

\bibliographystyle{agsm}

\bibliography{Clustertest}
\end{document}